\documentclass{article}

\usepackage{graphicx}
\title{\bf Numerical simulation of salt migration\\ -- Large deformation in viscoelastic\\ solid bodies}

 \author{\\ \\ I-Shih Liu$^1$%
 \footnote{Corresponding author: I-Shih Liu, e-mail: liu@im.ufrj.br} , \,
 Rolci A.\ Cipolatti$^1$,\\ Mauro A.\ Rincon$^1$, Luiz A.\ Palermo$^2$ %
 \\ \\
 $^1$Instituto de Matem\'{a}tica,
 Universidade Federal\\ do Rio de Janeiro, Brasil\\
 $^2$CENPES/Petrobras, Rio de Janeiro, Brasil}
 \date{}

\begin{document}
\maketitle

\begin{abstract}
We consider instability of a two layered solid body of a denser material on
top of a lighter one. This problem is widely known to geoscientist in
sediment-salt migration as salt diapirism. In the literature, this problem has
often been treated as Raleigh-Taylor instability in viscous fluids instead of
solid bodies. In this paper, we propose a successive linear approximation
method for large deformation in viscoelastic solids as a model for salt
migration.


\end{abstract}

 \newwrite\eqnfile
  \newif\ifeqnlist \newif\ifeqnwarning \newif\ifshoweqlabel
  \eqnlistfalse\showeqlabelfalse \eqnwarningtrue
\def\openeqnfile{\ifeqnlist\else\eqnlisttrue
    \immediate\openout\eqnfile=\jobname.eqn\fi}
\def\eqn(#1){\expandafter
     \ifx\csname#1eq\endcsname\relax
     \else\ifeqnwarning
    \immediate\write16{ * Equation (#1) is already defined.}\fi\fi
     \global\addtocounter{equation}{1}
     \expandafter\xdef\csname#1eq\endcsname{\theequation}
     \ifeqnlist\immediate\write\eqnfile{
    let (\csname#1eq\endcsname) = (#1).}\fi
     \eqprint{\theequation}
     \ifshoweqlabel\rlap{\scriptsize\hphantom{-}(#1)}\fi}%
\def\eqprint#1{(#1)}%
\def\eq(#1){\expandafter
     \ifx\csname#1eq\endcsname\relax
    \immediate\write16{ - Equation (#1) is undefined.}
    \eqprint{0.0}
     \else\eqprint{\csname#1eq\endcsname}\fi}%
\catcode `@=11
\def\eqalign#1{\null\,\vcenter{\openup\jot\m@th
  \ialign{\strut\hfil$\displaystyle{##}$&$\displaystyle{{}##}$\hfil
      \crcr#1\crcr}}\,}
\catcode `@=12
 \def\bfl#1{{\mbox{\boldmath$#1$}}}
 \def\bfs#1{{\setbox0=\hbox{$\scriptstyle#1$} 
     \kern-.015em\copy0\kern-\wd0
     \kern .030em\copy0\kern-\wd0
     \kern-.015em\raise.01em\box0 }}
 \def\grad{\mathop{\rm grad}}
 \def\div{\mathop{\rm div}}
 \def\Grad{\mathop{\rm Grad}}
 \def\Div{\mathop{\rm Div}}
 \def\tr{\mathop{\rm tr}}
 \def\ds{\displaystyle}
 \def\Re{I\!\!R}
 \def\Re{I\!\!R}
 \def\I{\hbox{$I$}}
 \def\II{\I\!\I}
 \def\Min{\mathop{\rm Min}\limits}
 \def\qed{{\mbox{{$\sqcup$}\llap{$\sqcap$}}}}
 \def\Line{\raise -.04ex\hbox{\bf---}}

\section{Introduction}
The behavior for large deformation is characterized by some nonlinear
constitutive relations, which leads to a system of nonlinear partial
differential equations. To solve boundary value problems of this nonlinear
system for large deformation, we propose a method in a successively updated
referential formulation of linear approximation.

The method is based on the well-known problem of small deformation superposed
on finite deformation in the literature \cite{GRS}. Roughly speaking, at each step, the constitutive functions are calculated at the present state of deformation which will be regarded as the reference configuration for the next state, and assuming the deformation to the next state is small, the constitutive functions can be linearized. In this manner, it becomes a linear problem from one state to the next state successively with small deformations.

Numerical simulations for elastic bodies with this method have been successfully compared with some exact solutions for large deformations \cite{LCR10,Liu-bend}. In this paper, we consider the instability of a two-layered solid body of a denser elastic material on top of a lighter viscoelastic one. This problem is widely known to geophysicists in sediment-salt migration as salt diapirism. Due to technical
difficulties in dealing with large deformations in solid bodies, the problem
has often been treated as Raleigh-Taylor instability in viscous fluids
(see for example, \cite{IZa,KSBM}) instead of solid bodies.

We will show in our numerical simulation with finite element method the
formation of diapirs with viscoelastic solid models. The mesh of the
calculation domain is updated as the body deforms at every step, however, up
to very large deformation we have found that re-meshing is almost unnecessary.

\subsubsection*{Remark 1:}
For numerical computation of large deformations in finite elasticity, to
obviate the difficulty of handling nonlinearities, incremental methods are
widely discussed, for example, in the books by Ciarlet \cite{Ciarlet}, Oden
\cite{Oden}, and Ogden \cite{Ogden}. These methods consist in letting the
boundary forces vary by small increments from zero to the given ones and to
compute corresponding approximate solutions by successive linearization. The
general idea seems to be of purely mathematical concern regarding
linearization between successive boundary value problems. The problems are
usually formulated with domains in the {\em initial} reference configuration, i.e., in (total) Lagrangian formulation. This is the essential difference from the method proposed in this paper with ``updated'' Lagrangian formulation, namely, the reference configuration is updated to the present state at each step. However, we refrain from calling the proposed method as an {\em updated Lagrangian formulation} (UL), because UL formulation is widely known in numerical methods albeit frequently with quite different tenets (see for example \cite{BLM}) in which the basic equations and numerical treatment of boundary conditions are quite different.

The present approach has an additional advantage that the assignment of
boundary data is simpler and straightforward. They are prescribed on the present state and the difference between prescribing the corresponding boundary data on the surface of the consecutive states is of higher order which is insignificant in the linear approximation. Therefore, unlike the usual incremental method in Lagrangian formulation, which needs to prescribe the incremental boundary data on the surface of the initial reference state with values depending on the deformed geometry of the body, we are totally free from any such complicated concerns (see \cite{LCR10,Liu-bend}).

\subsubsection*{Remark 2:}
For large deformations, systems of governing partial
different equations are generally nonlinear. The boundary value problems
are usually formulated in referential (Lagrangian) or spatial (Eulerian)  coordinates, and numerical solution of nonlinear systems of
algebraic equations by Newton's type methods. In this paper, we formulate boundary value problems in the coordinates relative to the configuration at the present time. In other words, we shall describe deformations relative to the current configuration, instead of a fixed reference configuration in the
Lagrangian formulation. Note that this is not an Eulerian formulation either.
The Eulerian coordinates are fixed coordinates in space. We shall call this a
relative-descriptional formulation \cite{Liu-bend}.

The advantage of using the relative-descriptional formulation is that we can
consider a small time step from the present state, so that the constitutive
functions can be linearized relative to the present state, and the equation of
motion becomes linear in relative displacement. When the present state
proceeds in time, a nonlinear finite deformation can be treated as a sequence
of small deformations in the same manner as the usual Euler's method for
solving differential equations, i.e., successively at each state, the tangent
is calculated and used to extrapolate the neighboring state.

The essential idea of the proposed method is based on the approach of small deformation superposed on large deformation \cite{GRS}. Although the {\it small-on-large} idea is very well-known, to my knowledge, numerical implementations of Euler's type successive approximation based on this idea have not been explored.

On the other hand, computational literature based on the {\it small-on-large} idea, typically the methods introduced in \cite{Simo}, \cite{SH}, and \cite{ZT} are very well documented. Such methods, in which the nonlinear algebraic systems resulting from the variational formulation of the {\it nonlinear} boundary value problem are solved by Newton's method, employ the tangent operator obtained from small-on-large linearization of the constitutive functions and boundary conditions. Such methods are quite different from the one proposed in this paper concerning the use of the small-on-large approach.

\section{The present configuration}

Let $\kappa_0$ be a reference configuration of the body ${\cal B}$, and
$\kappa_{t}$ be its deformed configuration at the {\em present} time $t$, Let
$$
  \bfl x = \chi(X,t)
  \qquad X\in\kappa_0({\cal B}),
$$
and
$$
  F(X, t)=\nabla_{X}^{}(\chi(X,t))
$$
be the deformation and the deformation gradient from $\kappa_0$ to
$\kappa_{t}$.

Now,  at some time~$\tau$, consider the deformed configuration
$\kappa_\tau$,
$$
  \bfl\xi=\chi(X,\tau):=\bfl\xi(\bfl x,\tau)\in\kappa_\tau{(\cal B)},
  \qquad \bfl x = \chi(X,t)\in\kappa_t{(\cal B)}.
$$
It can also be regarded as the relative deformation at time~$\tau$ with respect to the present configuration at time~$t$ denoted as $\bfl\xi(\bfl x,\tau)$. We also define the relative displacement vector as
$$
  \bfl u(\bfl x,\tau)=\bfl\xi(\bfl x,\tau)-\bfl x.
$$
Note that
$$
  \nabla_x\bfl\xi(\bfl x,\tau)= \nabla^{}_X(\chi(X,\tau))\,\nabla^{}_X(\chi(X,t))^{-1}
  =F(X,\tau)F(X,t)^{-1},
$$
hence, we have
$$
  \nabla_x\bfl u(\bfl x,\tau)=F(X,\tau)F(X,t)^{-1} - I,
$$
or simply as
$$
  F(\tau)=(I+H(\tau))F(t),
 \eqno\eqn(FF0)
$$
where $I$ is the identity tensor and
$$
 H(\bfl x,\tau) = \nabla_x\bfl u(\bfl x,\tau)
 \eqno\eqn(uH)
$$
is the relative displacement gradient at time $\tau$ with respect to the {\em present} configuration $\kappa_{t}$ (emphasize, not $\kappa_0$).
Moreover, by taking the time derivative with respect to~$\tau$, it gives
$$
  \dot F(\tau)=\dot H(\tau) F(t).
 \eqno\eqn(dotFF0)
$$
In these expressions and hereafter, we shall often denote a function~$F$ as $F(t)$ to emphasize its value at time~$t$ when its spatial variable is self-evident.

In summary, we can represent the deformation and deformation gradient
schematically in the following diagram:
\begin{center}
\begin{picture}(300,80)(-150,-40)

 \put(-30,-30){\vector(1,0){60}}
 \put(-10, 20){\vector(-1,-1){40}}
 \put( 10, 20){\vector( 1,-1){40}}

 \put(  0, 25){\makebox(0,10){$X\in\kappa_0({\cal B})$}}
 \put(-60,-35){\makebox(0,10){$\bfl x\in\kappa_t({\cal B})$}}
 \put(60,-35){\makebox(0,10){$\bfl\xi\in\kappa_\tau({\cal B})$}}

 \put(0,-25){\makebox(0,14)[b]{${I+H(\tau)}$}}
 \put(-45,0){\makebox(0,12)[t]{${F(t)}$}}
 \put(45,0){\makebox(0,12)[t]{${F(\tau)}$}}
 \put(0,-45){\makebox(0,14)[b]{$\xi=\bfl x+\bfl u$}}

\end{picture}
\end{center}

\section{Linearized constitutive equations}

Let $\kappa_0$ be the preferred reference configuration of a viscoelastic body
${\cal B}$, and let the Cauchy stress $T(X,t)$ be given by the constitutive
equation in the configuration $\kappa_0$,
$$
  T(X,t) = {\cal T}(F(X,t),\dot F(X,t)).
 \eqno\eqn(2.1)
$$
For large deformations, the constitutive function ${\cal T}$ is generally a
nonlinear function of the deformation gradient~$F$.

We shall regard the present configuration $\kappa_t$ as an {\em updated}
reference configuration, and consider a small deformation relative to the
present state $\kappa_t(\cal B)$ at time $\tau=t+\Delta t$. In other words, we
shall assume that the relative displacement gradient $H$ is small, $|H|\ll 1$, so that we can linearize the constitutive equation \eq(2.1) at time~$\tau$ relative to the updated reference configuration at time~$t$, namely,
$$\eqalign{
  T(\tau) &= {\cal T}(F(\tau), \dot F(\tau)) = {\cal T}(F(t),0) \cr
    &+ \partial_F{\cal T}(F(t),0)[F(\tau)-F(t)]+
     \partial_{\dot F}{\cal T}(F(t),0)[\dot F(\tau)]+o(2),}
$$
or by the use of \eq(dotFF0),
$$
    T(\tau) = T_e(t) + \partial_F{\cal T}(F(t),0)[H(\tau)F(t)] +
     \partial_{\dot F}{\cal T}(F(t),0)[\dot H(\tau) F(t)]+o(2).
$$
where
$$
  T_e(t)={\cal T}(F(t),0)
$$
is the elastic Cauchy stress at the present time~$t$ and $o(2)$ represents
higher order terms in the small displacement gradient $|H|$.

The linearized constitutive equation can now be written as
$$
  T(\tau) = T_e(t) + L(F(t))[H(\tau)] + M(F(t))[\dot H(\tau)],
 \eqno\eqn(T-cauchy)
$$
where
$$
  L(F)[H]: = \partial_F{\cal T}(F,0)[HF], \qquad
  M(F)[\dot H]: = \partial_{\dot F}{\cal T}(F,0)[\dot HF],
  \eqno\eqn(L-elast)
$$
define the fourth order {elasticity tensor} $L(F)$ and {viscosity tensor}
$M(F)$ relative to the present configuration $\kappa_{t}$.

The above general definition of the elasticity and the viscosity tensors for any constitutive class of viscoelastic materials $T={\cal T}(F,\dot F)$, relative to the updated present configuration, will be explicitly determined in the following sections for a particular class, namely a Mooney-Rivlin type materials.

\subsection{Compressible and nearly incompressible bodies}

For a viscoelastic body, without loss of generality, the constitutive equation
\eq(2.1) relative to the preferred reference configuration $\kappa_0$ can be
written as
$$
  T = {\cal T}(F,\dot F)
  = -p(F)I+\widetilde{\cal T}(F,\dot F).
 \eqno\eqn(2.1a)
$$
However, for an incompressible body, the pressure~$p$ depends also on the
boundary conditions, and because it cannot be determined from the deformation
of the body alone, it is called an {\em indeterminate pressure}, which is an
independent variable in addition to the displacement vector variable for
boundary value problems. \smallskip

For compressible bodies, we shall assume that the pressure depend on the
deformation gradient only through the determinant, or by the use of the mass
balance, depend only on the mass density,
$$
  p=\hat p(\det F)=p(\rho), \qquad
  \rho = {\rho_0\over \;\det F\;},
$$
where $\rho_0$ is the mass density in the reference configuration $\kappa_0$.
We have
$$\eqalign{
  \rho(\tau)-\rho(t)&=\rho_0(\det F(\tau)^{-1}-\det F(t)^{-1})
  =\rho(t)(\det (F(\tau)^{-1}F(t))-1) \cr
  &=\rho(t)(\det (I+H(\tau))^{-1}-1)
  =-\rho(t)\tr H(\tau) + o(2), \cr }
$$
in which the relation \eq(FF0) has been used.

Therefore, it follows that
$$
 p(\tau)-p(t) =\Bigl({dp\over d\rho}\Bigr)_{t}(\rho(\tau)-\rho(t))+o(2)
 =-\Bigl(\rho\,{dp\over d\rho}\Bigr)_{t} \,\tr H(\tau) +o(2),
$$
or
$$
  p(\tau)=p(t) -\beta\,\tr H(\tau) +o(2),
 \eqno\eqn(2.5)
$$
where $\beta:=\bigl(\rho\,{dp\over d\rho}\bigr)_{t}$ is a material parameter evaluated at the present time~$t$.
\medskip

From \eq(2.1a) and \eq(FF0), let
$$
  C(F(t))[H(\tau)] :=\partial_F^{}\widetilde{\cal T}(F(t),0)[F(\tau)-F(t)] =
  \partial_F^{}\widetilde{\cal T}(F(t),0)[H(\tau)F(t)],
$$
then from \eq(L-elast)$_1$, the elasticity tensor becomes
$$
  L(F)[H]=  \beta(\tr H)I + C(F)[H].
 \eqno\eqn(2.7)
$$

We call a body nearly incompressible if its density is nearly insensitive to
the change of pressure. Therefore, if we regard the density as a function of
pressure, $\rho=\rho(p)$, then its derivative with respect to the pressure is
nearly zero. In other words, for nearly incompressible bodies, we shall assume
that~$\beta$ is a material parameter much greater than~1,
$$
  \beta \gg 1.
 \eqno\eqn(2.6)
$$

Note that for compressible or nearly incompressible body, the elasticity
tensor does not contain the pressure explicitly. It is only a function of the
deformation gradient and the material parameter~$\beta$ at the present
time~$t$.

\subsection{A viscoelastic material model}

As an example, we shall consider a nearly incompressible viscoelastic material
with the following constitutive equation,
$$
  T = -p\,I+\widetilde{\cal T}(F,\dot F),
$$
$$\eqalign{
  \widetilde{\cal T}(F,\dot F) & =  s_1 B + s_2 B^{-1} \cr
   & + \lambda (\tr D) I +
  2\mu_1 D + \mu_2(DB+BD) +\mu_3(DB^{-1}+B^{-1}D),}
 \eqno\eqn(ex.1)
$$
where $B=FF^T$ is the left Cauchy-Green strain tensor and $D$ is the rate of
strain tensor. The material parameters $s_1$ through $\mu_3$ are assumed to be
constants. This material model will be referred to as a {\em Mooney-Rivlin type isotropic viscoelastic solid}, since, by assuming the relevant material parameters to be constant, the elastic part of the constitutive
equation represents the well-known Mooney-Rivlin isotropic elastic solid in
nonlinear elasticity (see \cite{Liubook,NLFTM}). Note that this constitutive
equation contains only linear terms in~$D$ (see also \cite{Hutter}). Therefore, it would
be appropriate for large deformations with small strain rate.

After taking the gradients of $\widetilde{\cal T}(F,\dot F)$ with respect
to~$F$ at $(F,0)$, we have
$$
  C(F)[H] = s_1(HB+BH^T)-s_2(B^{-1}H +H^T B^{-1}),
$$
which by \eq(2.7), gives the elasticity tensor
$$
  L(F)[H] =  -\beta(\tr H)I + s_1(HB+BH^T)-s_2(B^{-1}H
  +H^T B^{-1}),
 \eqno\eqn(L-MR)
$$
and with respect to~$\dot F$, we obtain
$$
  M(F)[\dot H] = \lambda (\tr\dot H)I +M_0(\dot H+\dot H^T)+ (\dot H+\dot H^T)M_0,
 \eqno\eqn(M-MR)
$$
where
$$
  M_0 := {1\over2}\,(\mu_1 I+ \mu_2 B+ \mu_3 B^{-1}).
$$

From \eq(T-cauchy), we have
$$
  T(\tau) = T_e(t) + L(F(t))[H(\tau)] + M(F(t))[\dot H(\tau)].
 \eqno\eqn(T-Cauchy)
$$
Furthermore, the (first) Piola-Kirchhoff stress tensor at time~$\tau$ relative
to the present configuration at time~$t$, denoted by $T_t(\tau)$, is given by
$$\eqalign{
  T_t(\tau)
  & = \det(I+H)\,T(\tau) (I+H)^{-T} \cr
  & = \det(I+H)\bigl(
  T_e(t) + L(F)[H] +M(F)[\dot H] \bigr) (I+H)^{-T}  \cr
  & = (I+\tr H)\bigl(
  T_e(t) + L(F)[H]+M(F)[\dot H] \bigr) (I-H^T) + o(2)  \cr
  & = T_e(t) + (\tr H)T_e(t) - T_e(t)H^T + L(F)[H] +M(F)[\dot H] +o(2).
  }
$$
We can write the linearized Piola-Kirchhoff stress as
$$\eqalign{
  T_t(\tau) &= T_e(t) + (\tr H(\tau))T_e(t) - T_e(t)H(\tau)^T \cr
   & \quad+ L(F(t))[H(\tau)] +M(F(t))[\dot H(\tau)]. }
$$
Note that when $\tau\to t$, $H\to 0$, and hence $T_t(\tau)\to T(t)$,
therefore, the Piola-Kirchhoff stress, becomes the Cauchy stress at the
present time~$t$.

We can also write
$$
  T_t(\tau) = T_e(t) +  K(F(t),T_e(t))[H(\tau)] +M(F(t))[\dot H(\tau)],
 \eqno\eqn(T-Piola)
$$
where the Piola-Kirchhoff elasticity tensor is defined as
$$
 K(F,T_e)[H] := (\tr H)T_e - T_eH^T + L(F)[H].
 \eqno\eqn(K-Piola)
$$

In numerical examples presented later, the material is assumed to be of
Mooney-Rivlin type and these relations will be used.

\section{Boundary value problem}

Let $\Omega=\kappa_t({\cal B}) \subset\Re^3$ be the region occupied by the
body at the present time~$t$, and $\partial\Omega=\Gamma_1
\cup\Gamma_2$, and $\bfl n_\kappa$ be the exterior unit normal to
$\partial\Omega$. We can now consider the updated referential formulation of
the boundary value problem.

At time~$\tau>t$, we shall consider the boundary value problem in the
Lagrangian formulation, with the present state at time~$t$ as the updated
reference configuration, given by
$$\left\{ \eqalign{
   - \div T_t(\tau) &= \rho(t)\bfl g(\tau),
   \;\quad\mbox{in~} \Omega , \cr
  T_t(\tau)\,\bfl n_\kappa &= \bfl f(\tau),
   \quad\mbox{on~} \Gamma_1, \cr
  \bfl u(\tau)\cdot\bfl n_\kappa &= 0,
   \quad\mbox{on~} \Gamma_2, \cr
  T_t(\tau)\,\bfl n_\kappa\times\bfl n_\kappa &=0,
   \quad\mbox{on~} \Gamma_2,
 }\right.
 \eqno\eqn(bvp)
$$
where $T_t(\bfl x,\tau)$ is the Piola-Kirchhoff stress at time~$\tau$ with
respect to the state at the present time~$t$. The divergence operator is
relative to the coordinate system $(\bfl x)$ in the present configuration. The
body is subjected to the body force $\bfl g(\bfl x,\tau)$ and the surface
traction $\bfl f(\bfl x,\tau)$.

From the definition \eq(uH), we also have the initial condition for the
displacement vector $\bfl u(\bfl x,\tau)$:
$$
 \bfl u(\bfl x,t)=0, \qquad \forall\,\bfl x\in \kappa_t({\cal B}).
 \eqno\eqn(ICu)
$$
Since $ T_t(\tau)\bfl n_\kappa\times\bfl n_\kappa=0 $ implies that the surface
traction $ T_t(\tau)\bfl n_\kappa$ is in the direction of the normal, the last
boundary condition in \eq(bvp) states that the tangential component of the
surface traction $T_t(\tau)\bfl n_\kappa$ vanishes on $\Gamma_2$. In other
words, the boundary~$\Gamma_2$ is a roller-supported boundary.

\subsection{Linearized boundary value problem}
We shall assume that at the present time~$t$, the deformation gradient $F$
with respect to the preferred reference configuration $\kappa_0$ and the
Cauchy stress $T$ are known, and that $\tau=t+\Delta t$ with small enough
$\Delta t$, then from \eq(T-Piola), the equilibrium equation in \eq(bvp) can
be written as
$$
   - \div(K(F(t),T_e(t))[\nabla \bfl u(\tau)] +M(F(t))[\nabla{\dot \bfl u}(\tau)])
   = \div T_e(t) +\rho(t)\bfl g(\tau),
 \eqno\eqn(pde)
$$
which is a linear partial differential equation for the displacement
vector~$\bfl u(\bfl x,\tau)$. The right hand side is a known function.

The idea of formulating the boundary value problem in the form \eq(bvp) for
small time increment is similar to the theory of small deformations superposed
on finite deformations (see \cite{GRS,NLFTM}).

\subsection{Variational formulation}
The boundary value problem \eq(bvp) can be formulated as a variational
problem. Let us consider the Sobolev space of vector valued functions on
$\Omega$,
$$
  H^1(\Omega)=\{\bfl v:\Omega\to \Re^3 \;\vert\; \bfl v,\,\nabla \bfl v\in L^2(\Omega)\}
$$
and the subspace
$$
  V =\{ \bfl v\in H^1(\Omega)\;\vert\; \bfl v\cdot \bfl n_\kappa=0
  ~{\rm on}~ \Gamma_2\}.
$$

Taking the inner product of the equation \eq(pde) with a vector $\bfl w\in V$
and integrating over the domain~$\Omega$, we obtain, after integration by
parts,
$$\eqalign{
 &\int_\Omega \Bigl(
  K[\nabla\bfl u(\tau)] + M[\nabla\dot{\bfl u}(\tau)]\Bigr)
  \cdot\nabla\bfl w \, dv
   \cr
 &\qquad\quad = \int_\Omega \rho(t)\bfl g(\tau)\cdot\bfl w\,dv
  - \int_\Omega T_e(t) \cdot\nabla\bfl w\,dv
  + \int_{\partial\Omega} T_t(\tau)\bfl n_\kappa\cdot\bfl w\,da,}
$$
where the relation \eq(T-Piola) is used and dot ($\,\cdot\,$) represents the
inner product between two vectors as well as the inner product between two
second order tensors, i.e., $A\cdot B=\tr (AB^T)$.

Since $\bfl w\in V$, $\bfl w\cdot\bfl n_\kappa=0$, by the boundary conditions
the surface integral on the right hand side becomes
$$
  \int_{\partial\Omega} T_t(\tau)\bfl n_\kappa\cdot\bfl w\,da =
  \int_{\Gamma_1}\bfl f(\tau)\cdot\bfl w\, da.
$$
Therefore, if we define the following bilinear forms
$$\eqalign{
  &{\cal L}(\bfl w, \bfl u(\tau)):=
  \int_\Omega K[\nabla\bfl u(\tau)]\cdot\nabla\bfl w\,dv, \cr
  &{\cal M}(\bfl w, \bfl u(\tau)):=
  \int_\Omega M[\nabla{\bfl u}(\tau]\cdot\nabla\bfl w\,dv,}
  \eqno\eqn(LM)
$$
and the linear form,
$$
  {\cal P}(\bfl w) :=
  \int_\Omega \rho(t)\bfl g(\tau)\cdot \bfl w\,dv +
  \int_{\Gamma_1}\bfl f(\tau)\cdot\bfl w\, da -
  \int_\Omega T_e(t)\cdot\nabla\bfl w\, dv,
  \eqno\eqn(LP)
$$
then the variational problem is to find the solution vector~$\bfl u(\tau)\in
V$ such that
$$
 {\cal L}(\bfl w, \bfl u(\tau)) + {\cal M}(\bfl w, \dot\bfl u(\tau)) =
 {\cal P}(\bfl w) \quad\forall\,\bfl w\in V.
 \eqno\eqn(vp)
$$

Note that from the initial condition \eq(ICu), $\bfl u(\bfl x,t)=0$, we can
approximate
$$
  \dot\bfl u(\tau) \approx {1\over \Delta t}\,\bfl u(\tau),
$$
and hence restate the variational problem as: Find the solution vector~$\bfl
u(\tau)\in V$ such that
$$
  {\cal K}(\bfl w, \bfl u(\tau)) =  {\cal P}(\bfl w) \quad\forall\,\bfl w\in V,
 \eqno\eqn(VP)
$$
where
$$
 {\cal K}(\bfl w, \bfl u):= {\cal L}(\bfl w, \bfl u) +
 {1\over\Delta t} {\cal M}(\bfl w, \bfl u).
$$

The variational problem depends on the elasticity and viscosity tensor at the updated present state. Mathematical analysis of requirements at the present state for existence and uniqueness of solution will be presented in a forthcoming paper. In general, non-existence or non-uniqueness may occurred if such requirements are not fulfilled.

For numerical solutions of the variational equation, finite element method
will be used.

\section{Successive linear approximation}

Recall the Euler method of solving differential equation, say $\dot y=f(t)$,
that for a discrete time axis, $\cdots < t_{n-1}<t_n<t_{n+1} <\cdots$, and
$y(t_n)=y_n$, the solution curve can be constructed by $y_{n+1}= y_n +
f(t_n)\Delta t$, where $f(t_n)$ is the tangent of the solution curve at $t_n$.
We can use a similar strategy for solving problems of large deformation, by
solving linear variational problem stated in \eq(VP).

We consider a discrete time axis, $t_0<\cdots <t_n<t_{n+1} <\cdots$ with small
enough time increment $\Delta t$. Let $\kappa_{t_n}$ be the configuration of
the body at the instant $t_n$ and
$$
   \bfl x_n = \chi(X, t_n) \in \kappa_{t_n}({\cal B})
   \quad\mbox{for}\quad
   X \in \kappa_{r}({\cal B}), \quad
   \kappa_r=\kappa_{t_0}.
$$
Let the elastic Cauchy stress $T_e(\bfl x_n,t_n)$ and the deformation gradient
$F(\bfl x_n,t_n)$ relative to the preferred configuration $\kappa_0$ at the
present time $t=t_n$ be known. The boundary value problem \eq(bvp) at the
instant $\tau=t_{n+1}$ with the body force $\bfl b(\bfl x_n, t_{n+1})$ and the
surface traction $\bfl f(\bfl x_n, t_{n+1})$ in updated referential
formulation with respect to the present configuration $\kappa_{t_n}$, can now
be solved numerically as a problem in linear elasticity for the relative
displacement field $\bfl u(\bfl x_n,t_{n+1})$ from the present state at
${t_n}$.

After solving the problem \eq(bvp) at $t_n$, the configuration
$\kappa_{t_{n+1}}$ can be regarded as the reference configuration at the
updated present time $t_{n+1}$ from the displacement field, i.e.,
$$
 \bfl x_{n+1}=\chi(X,t_{n+1})=\bfl x_n +\bfl u(\bfl x_n,t_{n+1}),
$$
while the deformation gradient
$$
 F(\bfl x_{n+1},t_{n+1})
 = \bigl(I+\nabla \bfl u(\bfl x_n,t_{n+1})\bigr) F(\bfl x_n,t_n)
$$
and the elastic Cauchy stress $T_e(\bfl x_{n+1},t_{n+1})$ can be calculated
from \eq(T-Cauchy) at $t_{n+1}$ so that the updated referential formulation of
the problem in the form \eq(bvp), with the body force $\bfl b(\bfl x_{n+1},
t_{n+2})$ and the surface traction $\bfl f(\bfl x_{n+1}, t_{n+2})$, can
proceed again from the updated referential configuration at $t_{n+1}$. This
numerical procedure will be referred to as the {\em successively updated
referential formulation of linear approximation}, or simply as the method of
{\em Successive Linear Approximation} (SLA). The updating process can be
represented in the following schematic diagram:
$$
\begin{array}{cccccccccc}
  && {\lower 12pt\hbox{\rm BVP}} &&    \\
  \kappa_r({\cal B})
  -\!\!{\rm fintely\atop deforms}\!\!
  \rightarrow
  & \displaystyle{\kappa_t=\kappa_{t_n}\atop
  \kappa_{t_n}({\cal B})}&
  \Line\!\!\! 
  \longrightarrow &
  \displaystyle{\kappa_\tau=\kappa_{t_{n+1}}\atop
  \kappa_{t_{n+1}}({\cal B})}& \\ \\
  && \mbox{update} & \Downarrow \\
  &&& & {\lower 12pt\hbox{\rm BVP}} \\
  &&& \displaystyle{\kappa_t=\kappa_{t_{n+1}}\atop
  \kappa_{t_{n+1}}({\cal B})}&
  \Line\!\!\! 
  \longrightarrow
  & \displaystyle{\kappa_\tau=\kappa_{t_{n+2}}\atop
  \kappa_{t_{n+2}}({\cal B})}& \\ \\
  &&& &\mbox{update} & \Downarrow \\ \\
  &&& && \displaystyle{{\kappa_t=\kappa_{t_{n+2}}\atop\cdots\cdots}} \\
\end{array}
$$
where BVP stands for the linear boundary value problem formulated in \eq(bvp) to be solved successively after each update as shown in the above scheme, and the process starts with $\kappa_t=\kappa_{t_0}=\kappa_r$ which is the initial reference configuration.

\subsection{Incremental loadings}
Note that the force term ${\cal P}(\bfl w)$ in the variational equation
\eq(VP) defined in \eq(LP) is, in fact, a small quantity of the order of the
time increment $\Delta t$. Therefore, the method of SLA can be regarded as an
incremental method. Here we shall emphasize the incremental features of the
boundary value problem \eq(bvp).

For convenience, we shall denote the time dependence on $t_n$ as subindex~$n$
or simply as~$n$. For example, we write the Cauchy stress $T(t_n)=T_n$ and the
elastic Cauchy stress $T_e(t_n)=T_e(n)$, and we have
$$
  T_n = {\cal T}(F_n,\dot F_n) ={\cal T}(F_n,0) +
  \partial_{\dot F}{\cal T}(F_n,0)[\dot F_n]
  =T_e(n) +\partial_{\dot F}{\cal T}(F_n,0)[\dot F_n]\,.
 \eqno\eqn(Tn)
$$

Now, we can rewrite the variational equation \eq(VP) as,
$$
  {\cal K}(\bfl w, \bfl u_{n+1}) =  {\cal P}(\bfl w,n) \quad\forall\,\bfl w\in V,
 \eqno\eqn(V.P)
$$
where, from \eq(LP), the force term is given by
$$
  {\cal P}(\bfl w,n)=
  \int_\Omega \rho_n\bfl g_{n+1}\cdot \bfl w\,dv +
  \int_{\Gamma_1}\bfl f_{n+1}\cdot\bfl w\, da -
  \int_\Omega T_e(n)\cdot\nabla\bfl w\, dv.
 \eqno\eqn(P1)
$$
Note that by \eq(Tn), it follows that
$$
 \int_\Omega T_e(n)\cdot\nabla\bfl w\, dv =
 \int_\Omega T_n\cdot\nabla\bfl w\, dv -
 \int_\Omega \partial_{\dot F}{\cal T}(F_n,0)[\dot F_n]\cdot\nabla\bfl w\, dv,
 \eqno\eqn(P2)
$$
and
$$
 \int_\Omega T_n\cdot\nabla\bfl w\, dv =
  \int_{\partial\Omega} T_n\bfl n_\kappa\cdot\bfl w\, dv
  -\int_\Omega \div T_n\cdot\bfl w\, dv.
$$
On the other hand, noting that $T_{t_n}(t_n)=T_n$ is the Cauchy stress at $t_n$, from \eq(bvp) we have the following boundary value problem
at the present time~$t_n$,
$$\left\{ \eqalign{
   - \div T_n &= \rho_n\bfl g_n,
   \;\quad\mbox{in~} \Omega , \cr
  T_n\bfl n_\kappa &= \bfl f_n,
   \quad\mbox{on~} \Gamma_1, \cr
  T_n\bfl n_\kappa\times\bfl n_\kappa &=0,
   \quad\mbox{on~} \Gamma_2,
 }\right.
$$
which leads to
$$
 \int_\Omega T_n\cdot\nabla\bfl w\, dv =
  \int_{\Gamma_1} \bfl f_n\cdot\bfl w\, dv
  +\int_\Omega \rho_n\bfl g_n\cdot\bfl w\, dv.
 \eqno\eqn(P3)
$$
Combing \eq(P1), \eq(P2), and \eq(P3), we obtain
$$
  {\cal P}(\bfl w,n)= I_1(\bfl w,n) +I_2(\bfl w,n) +I_3(\bfl w,n),
 \eqno\eqn(P4)
$$
which contains three types of loading for the variational equation \eq(V.P),
namely,
$$\eqalign{
  I_1(\bfl w,n) &=\int_\Omega \rho_n(\bfl g_{n+1}-\bfl g_n)\cdot \bfl w\,dv, \cr
  I_2(\bfl w,n) &=\int_{\Gamma_1}(\bfl f_{n+1}-\bfl f_n)\cdot\bfl w\, da, \cr
  I_3(\bfl w,n) &=\int_\Omega \partial_{\dot F}{\cal T}(F_n,0)[\dot F_n]\cdot
  \nabla\bfl w\,dv.}
  \eqno\eqn(Int)
$$
The integral~$I_1$ represents the incremental body force and~$I_2$ is the
incremental surface traction between time steps~$t_n$ and~$t_{n+1}$. The third one $I_3$ is due to the viscous effect of the material body.

\subsection{Incremental approximation for large deformations}
If we assume that the functions $\bfl g(\bfl x,t)$ and $\bfl f(\bfl x,t)$ are
smooth in~$t$, then their increments, $\bfl g_{n+1}-\bfl g_n$ and $\bfl f_{n+1}-\bfl f_n$, are of the order of $\Delta
t=t_{n+1}-t_n$, and since the variational equation is linear, the solution
vector $\bfl u(\bfl x,t_{n+1})$ is also of the same order. For elastic
material bodies, these are the two possible types of incremental loading.

Starting from an initial solution and applying the method of SLA with proper
loading conditions at each time step, the problem of large deformation can be
obtained. Numerical examples employing the SLA method for incremental surface
traction of pure shear and of gradually bending a rectangular block into a
circular section have been considered in \cite{LCR10,Liu-bend} for Mooney Rivlin
elastic materials.

For viscoelastic material bodies, there is a possible loading due to the
integral~$I_3$. To see this, we shall consider the case without surface
traction $\bfl f=0$ and time-independent body force $\bfl g\,$ so that
$I_1=I_2=0$, hence the integral~$I_3$ is the only possible incremental
loading.

To begin with, for $n=0$, assume that at the initial time $t_0$, we have a
static equilibrium solution so that we have the initial conditions: $F_0=I$
and $\dot F_0=0$. Therefore, from \eq(V.P), \eq(P4) and \eq(Int)$_3$, we have
$$
  {\cal K}(\bfl w, \bfl u_{1}) = I_3(\bfl w,0)=0, \qquad\forall\,\bfl w\in V,
$$
which implies that the solution vector $\bfl u_1(\bfl x)=0$. Consequently,
$\dot{\bfl u}_1={1\over\Delta t}\bfl u_1 =0$ and $\dot H_1=\nabla\dot{\bfl u}_1=0$. Then, from~\eq(dotFF0), $\dot F_1=\dot
H_1F_0=0$, which, in turn, leads to
$$
  {\cal K}(\bfl w, \bfl u_{2}) = I_3(\bfl w,1)=0, \qquad\forall\,\bfl w\in V,
$$
and implies that $\bfl u_2$ must also vanish and so on. In other words, if the
initial solution is an equilibrium solution, the solution remains valid for
all time. However, this conclusion may not be true since the initial solution
may not be a stable equilibrium solution in general.

Therefore, in order to study the stability, a small perturbation of the
initial solution is needed so that $\bfl u_1\not=0$, and hence $I_3(\bfl
w,1)\not=0$, to trigger the successive evolution of deformations. Two such examples are considered in the following sections.

\section{Salt migration}
As an example of large deformation, we shall consider a body consisting of two
different layers initially, with the mass density of the upper layer, the
overburden sediment, greater than that of the bottom layer, the rock salt.

In a similar situation for viscous fluids, the inversion of density leads to
the so-called Rayleigh-Taylor instability due to buoyancy effect of gravity.
In the numerical simulation by the use of SLA method, we shall present the
results confirming the existence of similar instability for viscoelastic solid
bodies.

\begin{figure}[ht]
\begin{center}
 \includegraphics[width=4in]{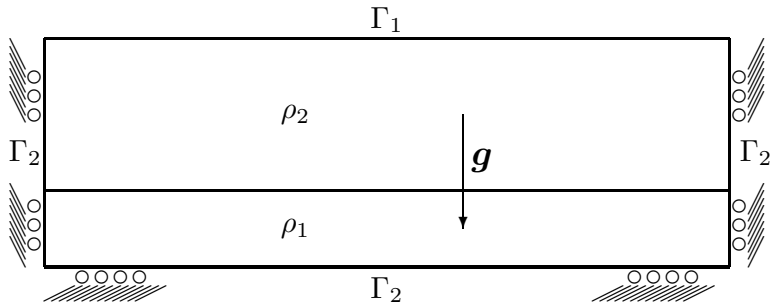}
\end{center}
 \caption{Boundary conditions and initial state.}
 \label{fig1}
\end{figure}

Consider a body consisting of two layers of elastic and viscoelastic solids as
shown in Fig.~\ref{fig1}.
The body is under the action of gravity $\bfl g$, and $\rho_1<\rho_2$. The upper boundary~$\Gamma_1$ is traction-free, the others~$\Gamma_2$ are
roller-supported. The initial state is an unstable equilibrium state, the
movement of the salt-sediment interface can be initiated by a small
perturbation.

For illustrative purpose, we shall present numerical simulations in a two-dimensional domain for a Mooney-Rivlin type material. The proposed method has been applied to three-dimensional domain and some different class of viscoelastic solid bodies with similar results.

\begin{figure}[p]
\begin{center}
 \includegraphics[width=4.7in,height=6.6in]{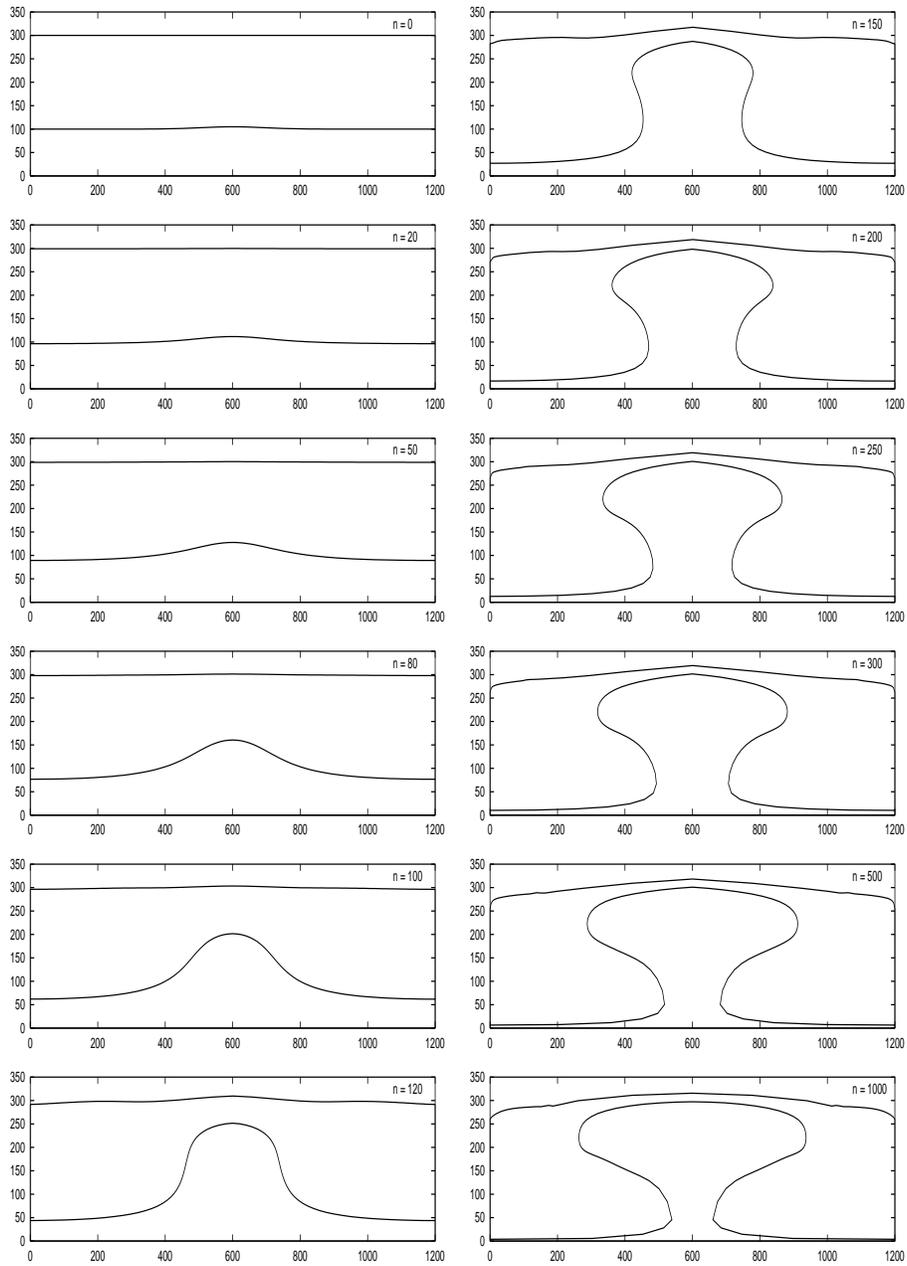}
\end{center}
 \caption{Migration of salt diapirism initiated by a small perturbation
 at the center of salt-sediment interface. The sequence represents the
 formation of a diapir during 100 million years ($t=n\,\Delta t$).}
 \label{fig2}
\end{figure}

\subsection{Formation of a salt diapir}\label{sec6.1}

In this example we consider formation of a salt diapir initiated by a small
perturbation at a very small region centered at the salt-sediment interface.
\begin{itemize}
\item
The dimension of the initial state is: length = 1,200 m, height of salt
layer~=~100~m, height of sediment layer = 200 m.
\item
The material parameters for rock salt are: \hfil\break
$\rho_0$ = 2.2$\times$10$^3$ Kg/m$^3$,  $s_1=0$, $s_2$ = -0.2$\times$10$^3$ Pa, \hfil\break $\lambda$~=~-10.0$\times$10$^3$~Pa\,Ma,
$\mu_1$~=~15.0$\times$10$^3$~Pa\,Ma, $\mu_2=\mu_3=0$, $\beta$~=~10$^9$~Pa.
\item
The material parameters for overburden sediment are: \hfil\break
$\rho_0$ = 3.0$\times$10$^3$ Kg/m$^3$,  $s_1$~=~2.5$\times$10$^3$~Pa,
$s_2$~=~-7.5$\times$10$^3$~Pa, \hfil\break
$\lambda=\mu_1=\mu_2=\mu_3=0$, $\beta$~=~10$^9$~Pa.
\item
The incremental time: $\Delta t$ = 0.1 Ma.
\end{itemize}
The material data are of convenient choice for demonstration only. No attempt
has been made to match the data to real properties of relevant materials.

\begin{figure}[h]
\begin{center}
 \includegraphics[width=4.7in,height=3.5in]{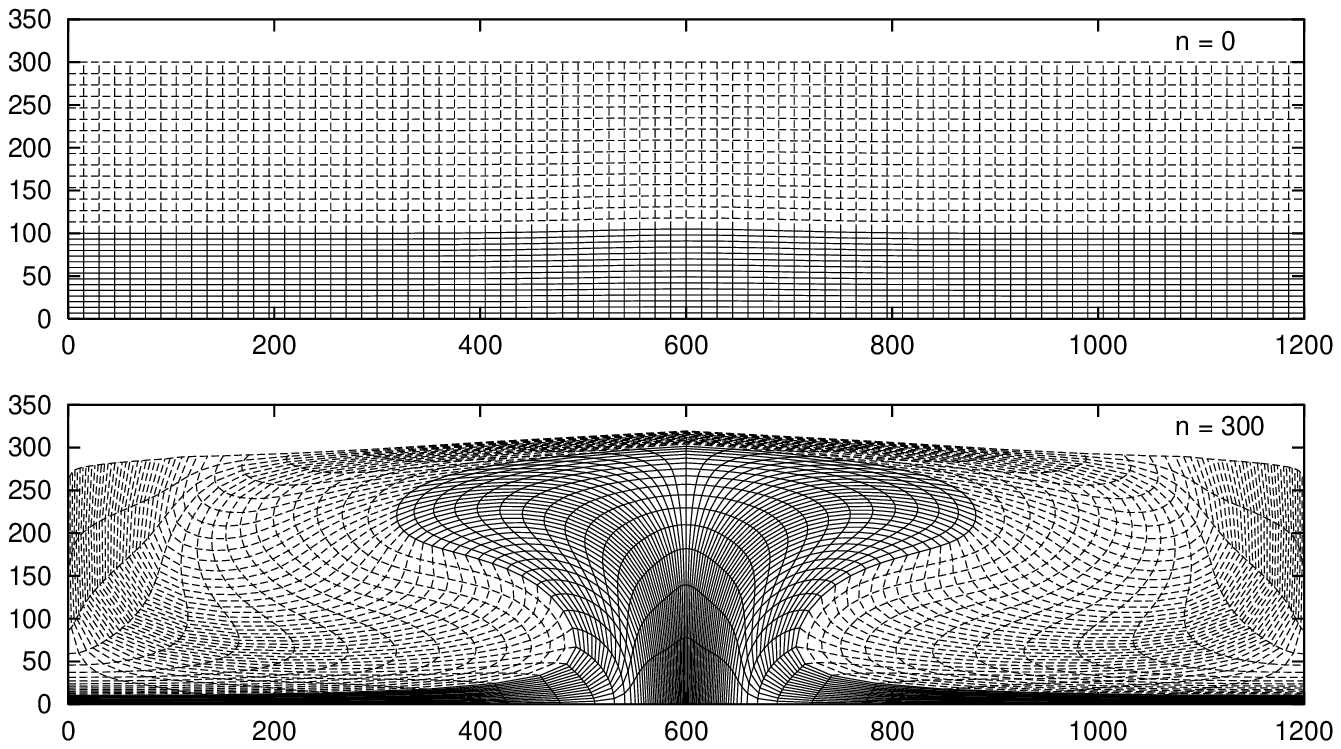}
\end{center}
 \caption{Original mesh and the deformed mesh at $n=300$.}
 \label{fig3}
\end{figure}

In Fig.~\ref{fig2}, various stages of migration of the salt-sediment interface
are shown, where~$n$ is the time step. One can easily see the formation of
salt diapir as time step increases. The effect is primarily due to buoyancy
force of density inversion. The diapir reaches its maximum height at about 20
Ma (million years) at $n$ = 200 with almost no change afterward as the diapir
becomes mature. The formation of diapir is a result of very large deformation
of the initial mesh as can be seen from Fig.~{\ref{fig3} at $n$ = 300. The
deformed mesh is quite similar to the experimental results of silicone putty
model of a diapir formed by spinning the model in a centrifuge by Dixon
\cite{Dixon}.

\begin{figure}[p]
\begin{center}
 \includegraphics[width=4.7in,height=6.6in]{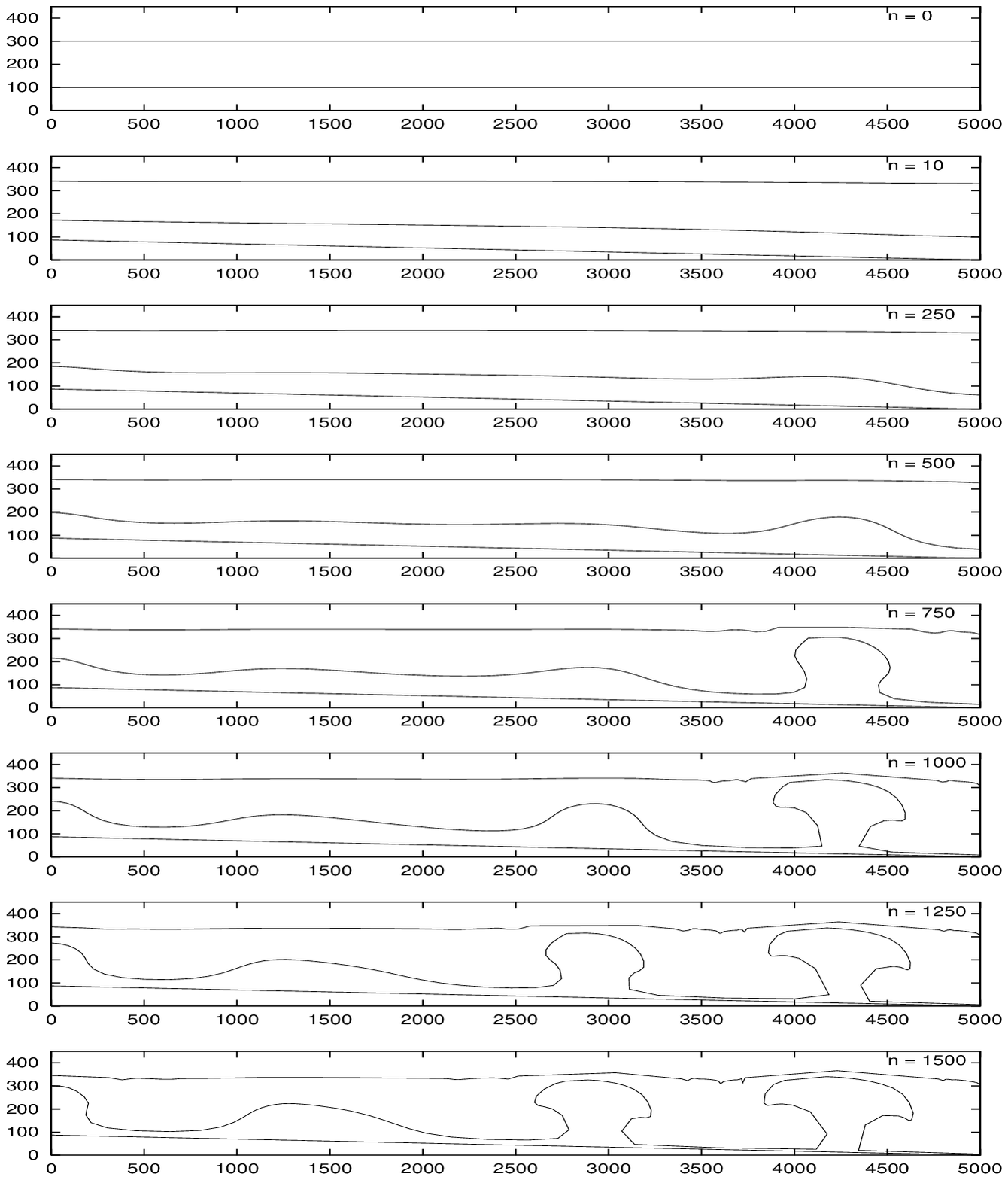}
\end{center}
 \caption{Salt migration initiated by a gradual lift of base rock to an angle
 of one degree at the initial 10 time steps on the left side. The figures
 show the growth sequence from right to the left for the appearance of salt
 structures.}
 \label{fig4}
\end{figure}

\subsection{Salt migration due to inclination}\label{sec6.2}

In the second example, we consider two-layer structure of a greater extension,
so that there are enough rock salt in the bottom layer to develop
multi-diapirs. The migration is initiated by gradually lifting the base rock,
which supports the salt layer, up to an inclination at an angle of one degree
during the initial one million years ($n$ = 10).
\begin{itemize}
\item
The dimension of the initial state is: length = 5,000 m, height of salt
layer~=~100~m, height of sediment layer = 200 m.
\item
The material parameters for rock salt are: \hfil\break
$\rho_0$ = 2.2$\times$10$^3$ Kg/m$^3$, $s_1=0$, $s_2$~=~{-0.2}$\times$10$^3$~Pa, \hfil\break
$\lambda$~=~-10.0$\times$10$^3$~Pa\,Ma,
$\mu_1$~=~15.0$\times$10$^3$~Pa\,Ma, $\mu_2=\mu_3=0$,
$\beta$~=~2$\times$10$^9$~Pa.
\item
The material parameters for overburden sediment are: \hfil\break
$\rho_0$ = 3.0$\times$10$^3$ Kg/m$^3$, $s_1$~=~2.5$\times$10$^3$~Pa,
$s_2$~=~-7.5$\times$10$^3$~Pa, \hfil\break
$\lambda=\mu_1=\mu_2=\mu_3=0$,
$\beta$~=~2$\times$10$^9$~Pa.
\item
The incremental time: $\Delta t$ = 0.1 Ma.
\end{itemize}

Due to the gravity, the lifting of the base rock on the left side pushes the
body to the right which initiates the growth of a salt pillow at about the
first 50 million years ($n=500$) as can be seen from Fig.~\ref{fig4}. As time
goes on, an adjacent salt pillow appears as the first one becomes a salt
diapir. The figures show the growth sequence from right to left of the
appearance of different salt structures up to 150 million years ($n=1500$).

\begin{figure}[h]
\begin{center}
 \vspace{0.3cm}
 \includegraphics[width=3.5in]{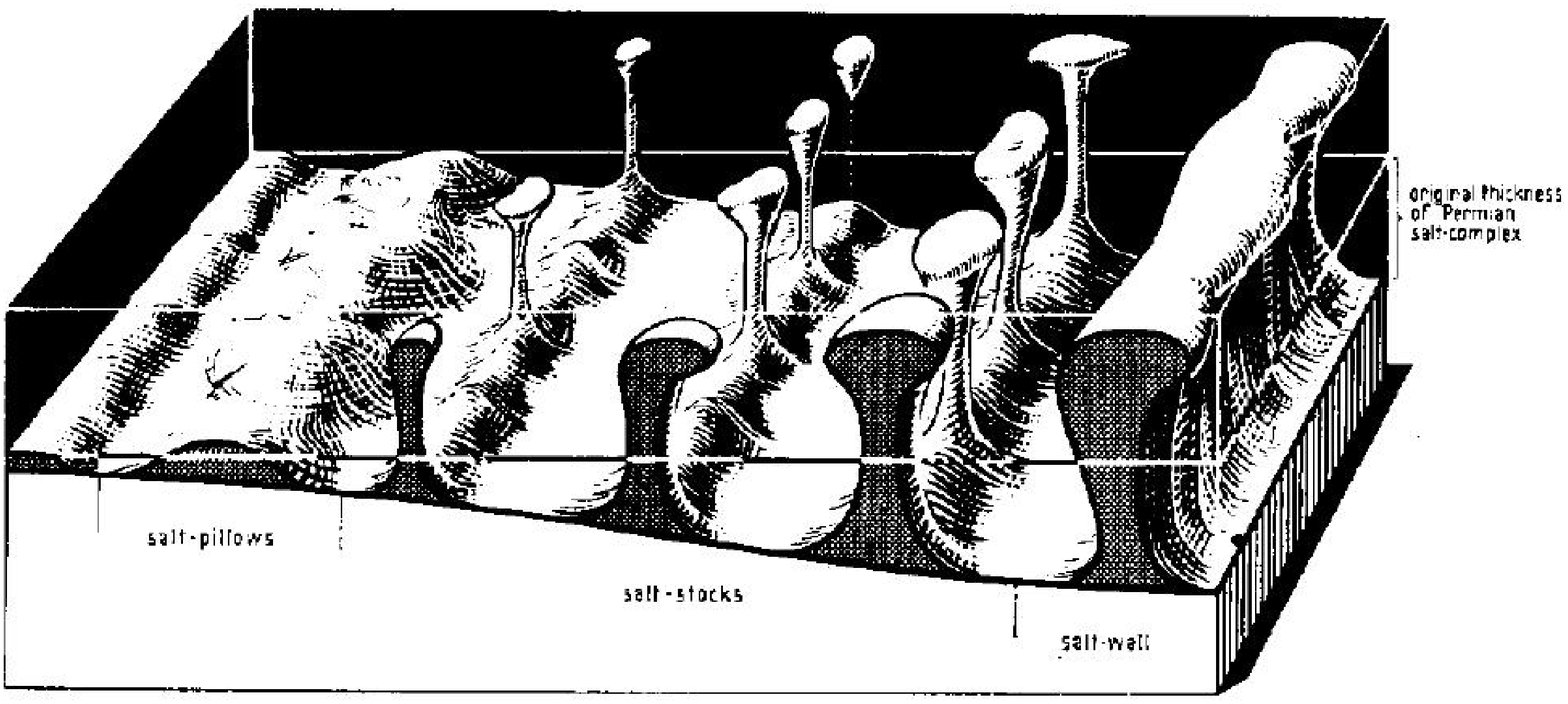} %
\end{center}
 \caption{Diagram of different types of salt structure due to a dip at an
 angle of more than one degree in Northern Germany (From
 F. Trusheim, Bull. Amer. Assoc. Petroleum Geologists 44 (1960))}
 \label{fig5}
\end{figure}

Such formation due to inclination of base rock was reported in permian salt
complex of northern Germany as shown in the sketch (Fig.~{\ref{fig5}}) by Trusheim \cite{Trusheim}. Although our numerical simulation is only two-dimensional, the similarity of the formation of salt structure is rather striking. However, since the numerical data are of convenient choice only, the estimated time in million of years in numerical simulation may not be of any real significance.

\medskip
\noindent{\bf Acknowledgement} {\small This work is supported by a research
project from {Petrobras/Brasil}. The authors (ISL, MAR) also acknowledge the
partial support from research fellowship of CNPq, Brazil.}

\end{document}